\documentclass[10pt]{amsart}

\usepackage[latin1]{inputenc}

\usepackage{graphicx,color}
\usepackage{wrapfig}
\usepackage{multirow}
\usepackage{array}
\usepackage{multicol}
\usepackage{amsmath,amssymb}
\usepackage{float}
\usepackage{mathrsfs}

\usepackage[all]{xy}
\xyoption{matrix}
\xyoption{arrow}

\def\ds{\displaystyle}

\def\complexes{\mathbb{C}}
\def\entiers{\mathbb{N}}
\def\naturels{\entiers}
\def\projectif{\mathbb{P}}

\def\reels{\mathbb{R}}
\def\sphere{\mathbb{S}}
\def\tore{\mathbb{T}}
\def\relatifs{\mathbb{Z}}

\def\ph{\varphi}

\def\rah{\begin{eqnarray}}
\def\arh{\end{eqnarray}}

\def\rae{\begin{eqnarray*}}
\def\are{\end{eqnarray*}}

\def\demi{\frac{1}{2}}
\def\ii{\sqrt{-1}}
\def\fr{\frac}
\def\dpe{\mathrel{\mathop:}=}
\def\tend{\rightarrow}
\def\sans{\backslash}

\def\commutatif{\ar@{}[rd]|{\circlearrowleft}}

\newcommand{\barre}[1]{\overline{#1}}
\newcommand{\inv}[1]{\frac{1}{#1}}

\DeclareMathAlphabet{\mathpzc}{OT1}{pzc}{m}{it}
\renewcommand{\Re}{\mathpzc{Re}\,}
\renewcommand{\Im}{\mathpzc{Im}\,}

\def\parel{\Re} 
\def\parim{\Im}

\newtheorem{thm}{\bf Theorem}[section]
\newtheorem{lem}[thm]{\bf Lemma}
\newtheorem{defn}[thm]{\bf Definition}

\newtheorem{prop}[thm]{\bf Proposition}

\newtheorem{cor}[thm]{\bf Corollary}

\newtheorem{rem}[thm]{\bf Remark}

\def\nin{\notin}

\usepackage{fullpage}

\def\preuve{\ \\\textit{Proof }}

\begin{document}

\title{Connectedness of the Arnold tongues for double standard maps}
\begin{author}[A.~Dezotti]{Alexandre Dezotti}
\email{dezotti$@$math.univ-toulouse.fr}
\address{ %
  Universit\'e de Toulouse; UPS; Institut de math\'ematiques de Toulouse UMR5219 \\ 
 118, route de Narbonne \\
  31062 Toulouse Cedex \\
  France }
\end{author}

\begin{abstract}
 We show that Arnold tongues for the family of double standard maps 
$$f_{a,b}(x)=2x+a-(b/\pi)sin(2 \pi x)$$
 are connected. This proof is accomplished in the complex domain by means of quasiconformal techniques and depends partly upon the fact that the complexification of $f_{a,b}$,
 has only one critical orbit taking symmetry into account.
\end{abstract}

\maketitle

\section{Introduction}

Our object is to study the family of mappings from the circle to itself,
 which we will call double standard maps,
 described below.
The starting point is the family of real analytic mappings
 $F_{a,b}:\reels\rightarrow\reels$
 given by :
\rae
X\mapsto F_{a,b}(X) =2X+a-\fr{b}{\pi}\sin(2\pi X),& & 
\are
 with $(a,b)\in\reels\times[0,1]$.
These functions descend to mappings 
 $f_{a,b}:\tore^1\rightarrow\tore^1$, where $\tore^1=\reels/\relatifs$.
We call this the double standard family
 $(f_{a,b})_{a\in\tore^1,b\in[0,1]}$.
This is a family of degree $2$
 mappings
 of the circle which has been recently studied by Ana Rodrigues and Micha\l \, Misiurewicz
 \cite{Misurod1},\cite{Misurod2}.

Each of $f_{a,b}$ in the family is monotonic semiconjugate to the doubling map
 $D:x\mapsto 2x \mod 1$
 via a unique mapping 
 $\ph_{a,b}$
(see lemma \ref{unicite de la semiconjuguaison}).
When we say a continuous mapping of the circle is monotonic we mean that its lift is a monotonic mapping of the real line.
This mapping
 $\phi_{a,b}$
 sends the periodic points of
 $f_{a,b}$
 on the periodic points of the doubling map
 $D$
 with same period.
Moreover
 $(a,b)\mapsto\phi_{a,b}$
 is continuous (see \cite{Misurod1}).

\ \\

In the following picture, graphs
 $(x,a(x))$
 of the parameters
 $(a,1)$
 for which
 $x$
 is a periodic point of
 $f_{a,1}$
 (with period $\leq 10$)
 are drawn in black.

\begin{center}
\includegraphics[scale=0.06]{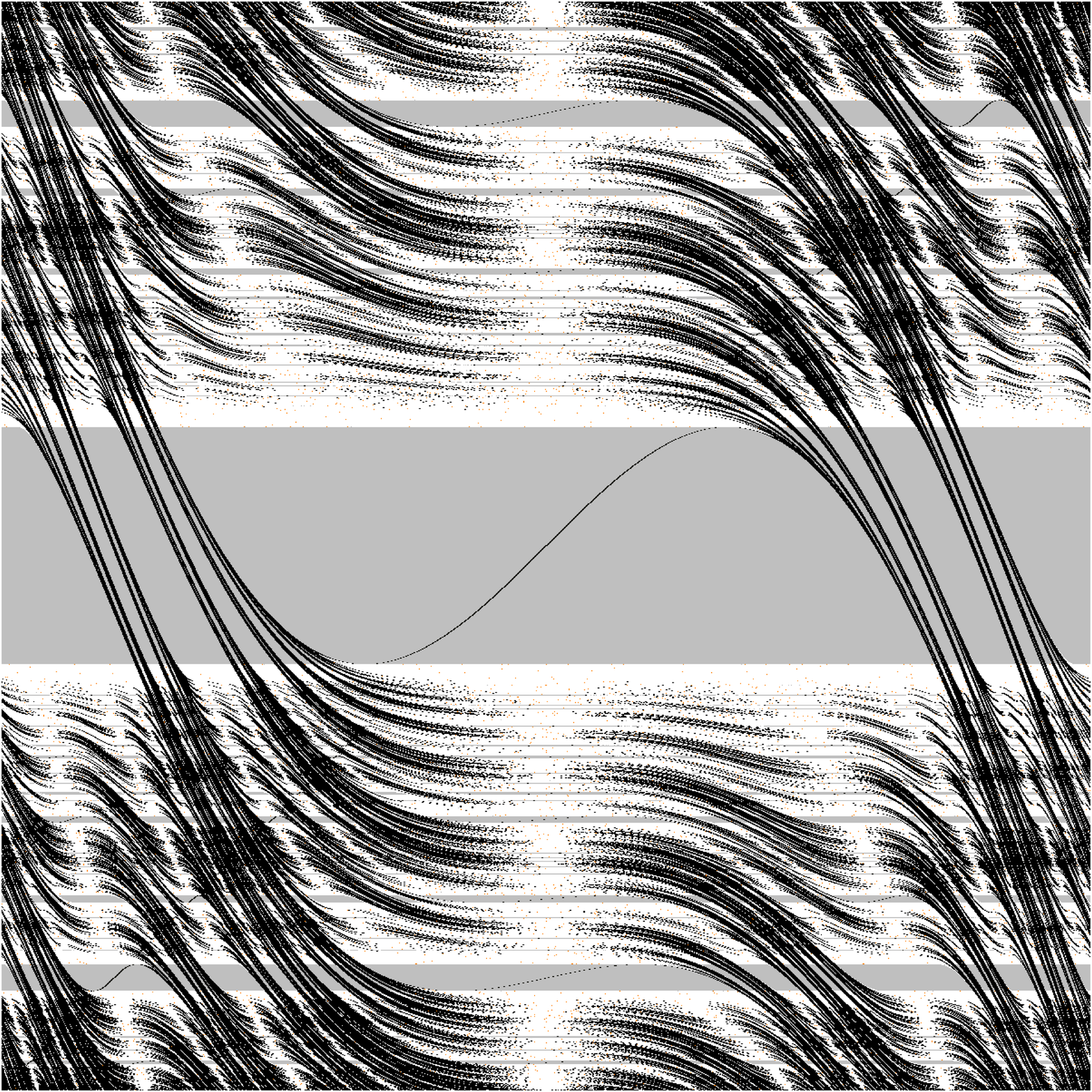} 
\end{center}
Indeed, all the functions
 $a(x)$
 are defined in
 $\tore^1$.
Part of the graphs, where the function decreases correspond to a repulsive cycle,
 whereas increasing parts correspond to an attracting cycle.
The intervals of values of
 $a$
 for which 
 $f_{a,1}$
 has an attracting cycle with period less than $10$ are emphasized in grey.

\ \\

One of the important fact concerning the mappings
 $f_{a,b}$
 is that they correspond to mappings\linebreak
 $g_{a,b}:\sphere^1\rightarrow \sphere^1$ of the complex circle (after conjugating by the exponential map) that
 extend to holomorphic functions
 $g_{a,b}:\complexes^*\rightarrow \complexes^*$
 given by
\rae 
g_{a,b}(z)=e^{2\ii\pi a}z^2e^{-b\left(z-\inv{z}\right)}.
\are
These functions
 $g_{a,b}$
 are symmetric with respect to the circle
 $\sphere^1=e^{2\pi\ii\tore^1}\subset\complexes^*$
 in the sense that
 $g_{a,b}\left(1/\barre{z}\right)=1/\barre{g_{a,b}(z)}$.
Any $g_{a,b}$
 has at most one attracting cycle on the circle and if there is one then all the critical points belong to the same component of the immediate basin.
Since this component contains only one point of the cycle, it yields a naturally distinguished point of the cycle.

\begin{defn}
A tuple of parameters $(a,b)\in\tore^1\times[0,1]$
 is of type $\tau$ if
 $g_{a,b}$ has an attracting cycle on the circle with
 $\phi_{a,b}(x_0)=\tau$, 
where $x_0$ 
 is the distinguished point of the attracting cycle of $g_{a,b}$.
\end{defn}
\begin{defn}
The tongue $T_\tau$ of type $\tau\in\tore^1$ is the set of parameters $(a,b)\in\tore^1\times[0,1]$
 of type $\tau$.
\end{defn}
\begin{rem}
Since any type $\tau$
 is a periodic point of $D$,
 it is a quotient
$\fr{k}{2^p-1}$ where $k\in\{0,\dots,2^p-2\}$
 and $p$ is the period of the cycle.
\end{rem}
The tongues
 $T_\tau$  are open disjoint subsets of
 $\tore^1\times[0,1]$.
 Some of their properties were studied in \cite{Misurod2}.

\ \\

The following picture shows the tongues (in black) in the parameter space
 $[-1/2,1/2]\times[0,1]$ :
\begin{center}
 \includegraphics[scale=0.18]{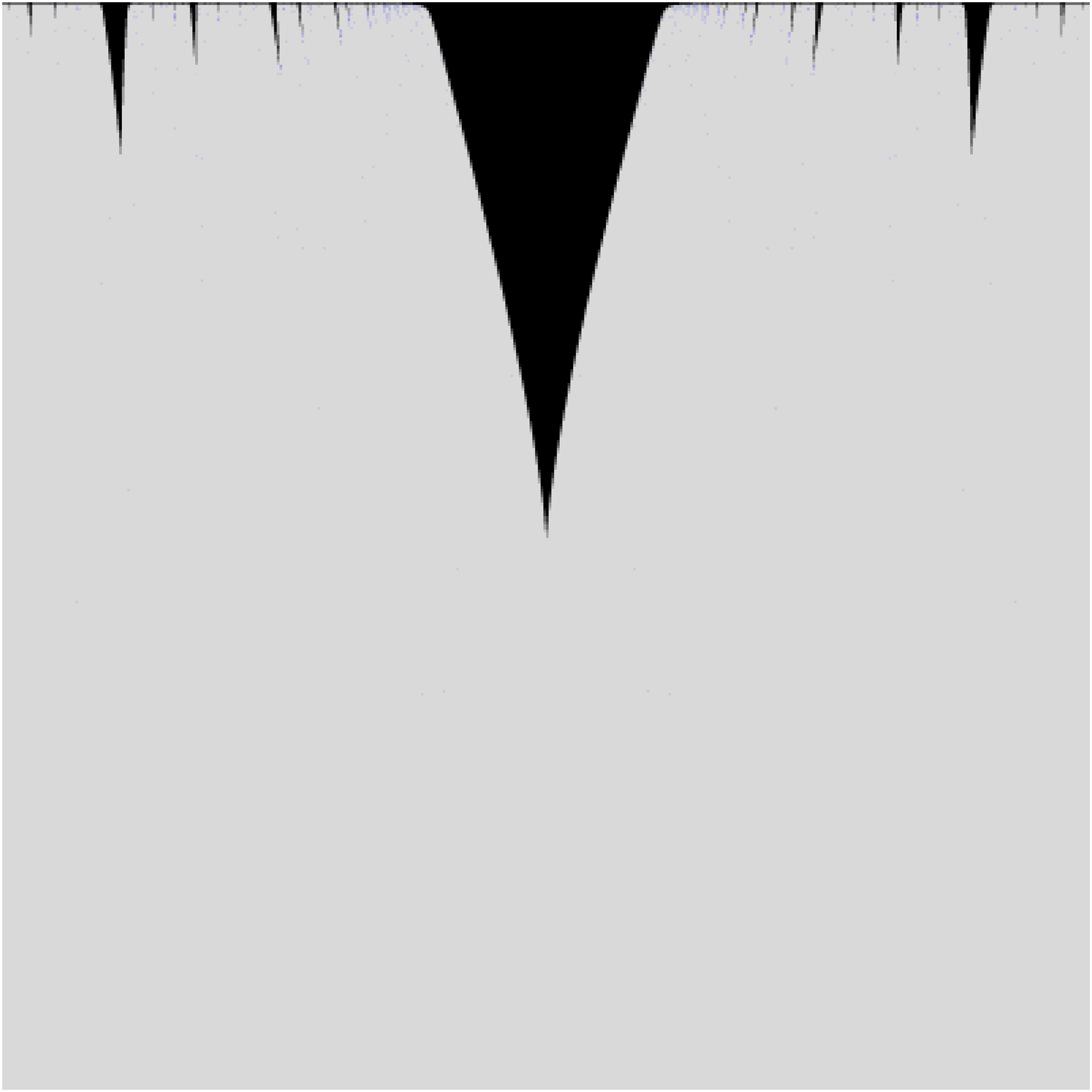} 
\end{center}

\ \\

One of the questions left unanswered by the authors of
 \cite{Misurod2} was about the connectedness of the set $T_\tau$.
We show here that the set
 $T_\tau$ is connected.
We shall use quasiconformal techniques in the proof.

\begin{thm}\label{thm chemin dans la langue et SA}
For each type $\tau$ (periodic point of $D$) there exists a unique 
 $(a_\tau,b_\tau)\in T_\tau$
 such that
 $f_{a_\tau,b_\tau}$ 
 has a superattracting cycle.
 Then 
 $b_\tau=1$,
 and all $(a,b)\in T_\tau$
 can be joined to $(a_\tau,1)$ by a continuous path in $T_\tau$.
\end{thm}

\begin{cor}
The tongues are connected.
\end{cor}

The paths will be constructed by a quasiconformal deformation.

Section \ref{section deformation}
 is devoted to the construction of the deformation inside the tongue for parameters with no superattracting cycle, see lemma \ref{chirurgie qc sur 2sm}.
In section \ref{section chemin dans la langue que c'est un chemain},
 it is shown that the deformation path tends to a unique limit parameter with a superattracting cycle.

\section{Acknowledgements}
The author wishes to thank Kuntal Banerjee, Xavier Buff, Adam Epstein, N\'uria Fagella and the referee.

\section{General properties of $g_{a,b}$}\label{sec gen prop of gab}

We recall several properties concerning of the mappings $g_{a,b}$ we use later on.

We firstly study the set of critical points of these mappings.
Let 
\rae
g(z) &=& z^2 e^{u(1/z)+v(z)},
\are
where $u$ and $v$ are polynomials of respective degrees $p\geq 1$ and $q\geq 1$.
Since $\deg v\geq 1$, the set of essential singularities of $g$ contains $0$.
The point $\infty$ is the only other essential singularity of $g$.

The derivative of $g$ is
$g'(z)  =  \fr{1}{z^{p-1}}P(z)e^{u(1/z)+v(z) }$,
where $P(z)=z^{p+1}v'(z)+2z^p-z^{p-1}u'(1/z)$.
The function $P(z)$ is a polynomial, its degree is $p+q$
 and its constant term is equal to the leading term of $u$,
 which is not zero. Therefore the critical points of $g$
 are the roots of $P(z)$ with same multiplicities.
The above fact is used in lemma \ref{lem le resultat est dans la famille important de tourner}.

Now if $p=q=1$, let $(a,b,c)\in \complexes/\relatifs\times \left(\complexes^*\right)^2$ be defined by
 $v(v)+u(1/z)=2\pi a\ii-(bz-c/z)$ ($a$ can only be defined modulo $1$).
 Then the critical points of $g$ are the roots of the quadratic polynomial
 $bz^2-2z+c$, namely $\fr{1\pm\sqrt{1-bc}}{b}$, where $\sqrt{1-bc}$
 stands for an arbitrary choice of one of the two square roots of $1-bc$.

We will say that a mapping $f:\complexes^*\rightarrow\complexes^*$
 is symmetric with respect to the circle if for all $z\in\complexes^*$,
$f\left(\inv{\barre{z}}\right)=\inv{\barre{f(z)}}$.
If we suppose that $g$ has the above symmetry property, then $c=\barre{b}$ and $a\in\reels$.

We now suppose that $g$ is symmetric with respect to the circle and denote it $g_{a,b}$, i.e.
 $g_{a,b}(z)=\lambda z^2 e^{-(bz-\barre{b}/z)}$, with $|\lambda|=1$.
One must notice that there exists a unique rotation that conjugates $g_{a,b}$
 to a function $g_{\tilde{a},\tilde{b}}=\tilde{\lambda}z^2 e^{-(\tilde{b}z-\tilde{b}/z)}$
 with $\tilde{b}\in\reels_+$.
Indeed, let $\rho$ be the unique complex number of modulus $1$ such that $\rho b\in\reels_+$,
 then 
 $\inv{\rho}g_{a,b}(\rho z)=\rho\lambda z^2 e^{-\left(\rho b z-\barre{\rho b}/z\right)}$ and
 $\tilde{b}=\rho b=\barre{\rho b}\in\reels_+$.
 Unicity follows from the above formula and from the unicity of $\rho$.

When $|b|<1$, the critical points of  $g_{a,b}$,
 which are $\left\{\fr{1-\sqrt{1-|b|^2}}{b},\fr{1+\sqrt{1-|b|^2}}{b}\right\}$,
 don't belong to the circle, from what follows that
 $g_{a,b|\sphere^1}$ is a monotonic mapping of the circle.
From the continuity dependance of $g_{a,b}$ on $b$, it follows that for any $a$, the mapping $g_{a,b|\sphere^1}$
 is also monotonic when $|b|=1$.

Conversely, suppose $|b|>1$,
 then $g_{a,b|\sphere^1}$
 is not monotonic. Here it is convenient to make the computations with the lifted map $f_{a,|b|}$
 of the mapping
 $g_{a,|b|}$, which is conjugated with $g_{a,b}$ by a rotation. We have
 $f_{a,|b|}(x)=2x+a+\fr{|b|}{\pi}\sin 2\pi x$. The critical points of $f_{a,|b|}$ are the solutions of
 $\cos(2\pi x)=\inv{|b|}$,
 these critical points all belong to the real axis and the second derivative does not vanish at these points.
Consequently $f_{a,|b|}$ and hence $g_{a,b|\sphere^1}$ are not monotonic when $|b|>1$.

We now compute the asymptotic values of the mapping $g_{a,b}$. Recall that $\alpha$
 is an asymptotic value of
 $g_{a,b}$ if there exists a continuous path $\eta:[0,+\infty[\rightarrow \complexes^*$ such that the limit of
 $\eta(t)$ when $t\tend +\infty$ is either $0$ or $\infty$ and such that
 $\alpha=\ds{\lim_{t\tend+\infty}}g_{a,b}(\eta(t))$.
\begin{prop}
Let $\lambda\in \sphere^1$, $b\in\reels^*$. Then the set of asymptotic values of
 $g(z)=\lambda z^2 e^{-b\left(z-1/z\right)}$ is $\{0,\infty\}$.
\end{prop}

\preuve
One easily check that $0$ and $\infty$ are asymptotic values.
We show there is no other asymptotic value. We use an argument inspired by \cite{ZakeriSiegel}.

Since we have $g(1/\barre{z})=1/\barre{g(z)}$, 
 we can consider only paths tending to $\infty$.

Let $\eta:[0,+\infty[\rightarrow \complexes^*$ be such a path and let $x\in\reels,y\in\reels$
 be such that $x+iy=-b\left(\eta(t)-1/\eta(t)\right)$.
Since we don't want $\infty$
 as  a limit point of $(g(\eta(t)))_{t\tend+\infty}$, we have to suppose that
 $x\tend -\infty$ as $t\tend +\infty$.
From the fact that $x\sim\parel(-b\eta)$ when $\eta\tend\infty$, it follows that
 $\eta$ cannot wind infinitively many times round $0$.
As a consequence we can define a continuous mapping $t\mapsto\log\eta(t)$
 with bounded imaginary part.

On the other hand, we must assume
 $y=\parim\left(-b\left(\eta(t)-1/\eta(t)\right)\right)$
 is unbounded. Indeed $\log|g|=2\log|\eta|+x$ is comparable to $\log|-b\eta|+x$, so if $y$ is bounded then
 $\log|g|$ is comparable to $x+\log|x|$, which tends to $-\infty$.
 This implies that $g\circ\eta\tend 0$ which we exclude.

From this, it follows that $\parim\left(\log g(\eta)\right)=\parim\left(2\log\eta\right)+y$ is unbounded and
 that the only possible limits of $g(\eta(t))$ when $t\tend+\infty$ are $0$ and $\infty$ $\square$

One could also use the Ahlfors-Carleman-Denjoy theorem to show that $g_{a,b}$
 has no asymptotic value in $\complexes^*$. The result follows from the fact that the entire function
 $f_{a,b}(z)=2z+a+\fr{b}{\pi}\sin 2\pi z$ has no finite asymptotic values.
Note that we would use the fact that $f_{a,b}$
 has finite growth order.
 Since we will need the definition of the growth order for the characterization of the double standard family, we recall it there.

\ \\

\begin{defn}
Let $a$ be a point of the Riemann sphere $\mathbb{P}^1$, $U$ a neighbourhood of $a$
 and $f:U\sans\{a\}\rightarrow \complexes$
 a holomorphic function defined in the punctured neighbourhood $U\sans\{a\}$
 of $a$
 such that it has an essential singularity at $a$.
Then the growth order of $f$ near $a$ is defined by
\rae
\limsup_{r\tend 0} -\fr{\log\log M_a(f,r)}{\log r},
\are
where $M_a(f,r)=\ds{\sup_{d(z,a)=r}|f(z)|}$ and $d$
 is the spherical distance.
\end{defn}
This definition is independent of change of chart in the domain.
For the family of mappings $(g_{a,b})_{a\in\reels,b\in\reels^*}$, the growth orders near $0$ and near $\infty$
 are both $1$.

\ \\

Now we return to the problem of finding attracting basins for $g_{a,b}$.

\begin{prop}
For any $(a,b)\in \tore^1\times[0,1]$,
the mapping $g_{a,b}$
 can have at most two attracting cycles, and if there are, both are of same period.
Moreover the immediate basin of an attracing cycle of
 $g_{a,b}$ contains at least one critical point.

If an attracting cycle of $g_{a,b}$
 belongs to the circle, it attracts both critical points and there is no other attracting cycle. In that case, the critical points belong to the same component of the immediate basin.
\end{prop}
\preuve
Here the singularites of $g_{a,b}$
 and of its inverse function are the same, as a consequence,
for all $n\in\naturels^*$, $g_{a,b}^{\circ n}$
 has the same asymptotic values as $g_{a,b}$.

It also implies that
 the immediate basin of an attracting cycle of 
 $g_{a,b}$
 contains at least one critical point.
Indeed, suppose the immediate basin does not contain any critical point.
Since the asymtpotic values don't belong to the domain of $g_{a,b}$,
 the boundary of the image $U$
 of the maximal extension of the inverse $\psi:D_R\rightarrow U$
 of the linearizing map corresponding to this attracting cycle contains either $0$ or $\infty$.
Moreover, from the semiconjugacy it follows that the boundary of
 $g_{a,b}(U)$ is a smooth Jordan curve contained in $U$.
The derivative of $g_{a,b}^{\circ n}$, where $n$ is the period of the cycle, is non zero on
 $\partial U\sans\{0,\infty\}$.
Thus $\partial U$
 is a union of one or two pieces of smooth curves with at least one of the points
 $0$ or $\infty$. The mapping $g_{a,b}^{\circ n}$
 is injective on these curves, which images are contained in the smooth Jordan curve
 $\partial g_{a,b}^{\circ n}(U)$.
 Hence, at least one of the end points of the images of these curves is an asymptotic value.
 This is impossible because  we have $\partial{g_{a,b}^{\circ n}(U)}\subset U$.

Since there are only two critical points, there are at most two attracting cycles.
Symmetry implies they both have the same period.
Moreover, if this cycle belongs to the circle, from the symmetry property of $g_{a,b}$,
 it follows that both critical points belong to its immediate basin.$\square$

\section{Deformation}
\label{section deformation}

\begin{lem}
\label{chirurgie qc sur 2sm}
Let
 $(a,b)\in\tore^1\times[0,1]$ 
 be such that
 $g_{a,b}$
 has an attracting cycle on the circle of type
 $\tau$
 with multiplier
 $\lambda$
 and take
 $\rho \in ]0,1[$.
Then there exists
 $\ph_{\rho}:\mathbb{P}^1\rightarrow \mathbb{P}^1$
  a quasiconformal homeomorphism of
 $\mathbb{P}^1=\bf{P}^1(\complexes)$
 fixing
 $0$
 and $\infty$,
 depending analytically on the real variable $\rho$, such that  $\ph_\lambda=id_{\mathbb{P}^1}$ and
 $g_{a,b}$ is conjugated with $g_{a_\rho,b_\rho}$
 by
 $\ph_{\rho}$
 where $g_{a_\rho,b_\rho}$
 has an attracting cycle on the circle of type
 $\tau$ with multiplier
 $\rho$ :
\begin{equation}
\xymatrix{
\complexes^* \ar[r]^{\ph_\rho} \ar[d]_{g_{a,b}}  & \complexes^* \ar[d]^{g_{a_\rho, b_\rho}} \\
\complexes^* \ar[r]_{\ph_\rho} &\complexes^*}
\end{equation}

\end{lem} 

Let us show lemma \ref{chirurgie qc sur 2sm}.
The aim is to construct a mapping
 $g_{\tilde{a},\tilde{b}}$,
 belonging to the complex double standard family, from
 $g:=g_{a,b}$,
 such that $g_{\tilde{a},\tilde{b}}$ has same type as $g:=g_{a,b}$
 but with a given multiplier
 $\rho\in]0,1[$. This is done by quasiconformal deformation.

\subsection{Deformation pocess}
\label{sec def process}

Let $\ph$ be a local linearizing map of $g^{\circ p}$,
 where $p$ is the period of the attracting cycle, defined in a neighbourhood of
  $x\in \sphere^1$, where $x$
 is the point of the cycle which lies in the same Fatou component as the critical points of
  $g$. 

We will suppose that the derivative of $\ph$ at $x$
 satisfies $|\ph'(x)|=1$, that it is tangent to the unit circle $\sphere^1$
 and points towards the anticlockwise direction. The mapping
 $\ph$ is uniquely determined by this normalization.
We have $\ph\circ g^{\circ p} = \lambda\ph$, $\ph(x)=0$ and
  $\ph'(x)=\ii x$.

It is easy to see that
 $\ph$
 is symmetric from what follows that if $D_R$
 is a disk centered in $0=\ph(x)$ on which $\ph^{ -1}$ is well-defined, then
 $U=\ph^{-1}(D_R)$ is stable by $z\mapsto 1/\barre{z}$
(cf. lemma \ref{La linearisante est symetrique elle aussi} below).

\begin{lem}
\label{La linearisante est symetrique elle aussi}
Suppose $f$
 is a holomorphic function defined in a neighbourhood of the circle, which is symmetric with respect to the circle, i.e. such that
 $f(1/\barre{z}) = 1/\barre{f(z)}$, and that it has an attracting cycle
 $(x,f(x),\dots,f^{\circ p-1}(x))$ 
 in $\sphere^1$ with multiplier $\lambda$.

Then $\lambda\in\reels$ and the linearizing mapping $\ph$ defined in a neighbourhood of $x$  
 and with $\ph'(x)=\ii x$
 staisfies $\barre{\ph(1/\barre{z})}=\ph(z)$
 and the preimages by $\ph$
 of any disk centered in $0$ are symmetric with respect to the circle.
\end{lem}

\preuve
Since $f$ is symmetric with respect to $\sphere^1$, the multiplier $\lambda$ is real.
 If we set $\psi(z)=\barre{\ph(1/\barre{z})}$ then it's easy to see that
 $\psi\circ f =\lambda\psi$ because $\lambda$ is real.
Moreover $\psi$ is holomorphic and a simple computation shows that
 $\psi'(x)=-\left(\inv{x}\right)^2\barre{\ph'(1/\barre{x})}=\ph'(x)$.

Therefore, from the uniqueness of $\ph$, it follows that $\ph(1/\barre{z}) = \barre{\ph(z)}$.

Then $\ph(z)\in D_R \Leftrightarrow \barre{\ph(z)}\in D_R\Leftrightarrow \ph( 1/\barre{z})\in D_R$.
$ \Box $

\ \\

We firstly use a quasiconformal deformation $\chi$ which transforms $D_R$ into
 $D_r$,
 with $R$
 chosen so as to produce a deformed map with the expected multiplier
 $\rho$.

\begin{lem}
\label{lem cycle et multiplicateur apres changemt}
\label{lem sur le changement de structure c du disque}
Let 
\begin{equation}
\begin{array}{crcl}
\chi : & \complexes^* & \longrightarrow & \complexes^* \\
        & z                  & \longmapsto     & |z|^\alpha z.
\end{array}
\end{equation}
Then :
\begin{itemize}
\item we have
\begin{equation}\label{muqui}
 \mu_{\chi} \dpe \fr{\partial\chi/\partial\barre{z}}{\partial\chi/\partial z} = \fr{\alpha/2}{1+\alpha/2}\fr{z}{\barre{z}}.
\end{equation}
In particular $|\mu_\chi|$ is constant and $||\mu_\chi||_\infty = \left|\fr{\alpha/2}{1+\alpha/2}\right|$;
\item $\chi$ is an invertible, multiplicative and radial map, i.e. 
 $\chi(re^{i\theta}) = \phi(r)e^{i\theta}$, where $\phi(r)=r^{\alpha+1}$ ;
\item if $\alpha=\fr{\log r}{\log R}-1$, $\chi$ sends the disk
 $D_R$ (of radius $R$) onto the disk $D_r$ ;
\item if moreover $R,r<1$, then
\begin{equation}
|\mu_\chi| = \fr{\left|1-\log r/\log R\right|}{1+\log r/\log R} <1.
\end{equation}
We have the following commutative diagram :
\begin{equation}
\xymatrix{
D_R \ar[d]_{z\mapsto \lambda z} \ar[r]_\chi  & D_r   \ar[d]^{z\mapsto \rho z} \\
D_R \ar[r]^\chi & D_r}
\end{equation}
with $\rho=\chi(\lambda)=\lambda^{1+\alpha}$ (in particular $0<\rho<1$).
\end{itemize}
\end{lem}

\preuve

\begin{itemize}
\item
 We have $\fr{\partial|z|}{\partial z} = \demi\fr{\barre{z}}{|z|}$ and 
 $\fr{\partial|z|}{\partial \barre{z}} = \demi\fr{z}{|z|}$.
Hence
$\fr{\partial\chi}{\partial\barre{z}} = \fr{\alpha}{2} |z|^{\alpha-2} z^2$ and
$\fr{\partial\chi}{\partial z} = \left(\fr{\alpha}{2}+1\right)|z|^\alpha$;


\item
the mapping $\chi$ is radial, monotonic on the radii, $\chi(0)=0$ and 
 $\chi(R)=R^{\fr{\log r}{\log R}} =r$;

\item $\chi(\lambda z) = \chi(\lambda) \chi(z)=\rho \chi(z)$ since $\chi$ is multiplicative.
$ \Box $
\end{itemize}
\ \\

\begin{lem}\label{lem sur le changement de structure c du disque num2}
The composition 
\begin{equation*}
U \stackrel{\ph}{\longrightarrow} D_R \stackrel{\chi}{\longrightarrow} D_r
\end{equation*}
produces a quasiconformal structure on 
 $U$, depending analytically on $\rho$, which is invariant by
 $g$ and
 with complex dilatation coefficient :
\begin{equation}
\mu_{\rho}  = 
\fr{\log r/\log R-1}{1+\log r/\log R}\cdot
\fr{\ph(z)}{\barre {\ph(z)}}\cdot\fr{\barre{\ph'(z)}}{\ph'(z)}
 =\fr{\alpha/2}{1+\alpha/2}\cdot\fr{\ph(z)}{\barre{\ph(z)}}\cdot\fr{\barre{\ph'(z)}}{\ph'(z)}
\end{equation}
satisfying $\mu_{\rho}(1/\barre{z})=\barre{\mu_{\rho}(z)}$ on  $U$.
\end{lem}
\ \\
\preuve
The complex dilatation $\mu_{rho}$ is the coefficient of the Beltrami form
 $\fr{\partial \chi\circ \ph}{\barre{\partial}\chi\circ\ph}$.
Analytic dependence follows from \ref{muqui} and the fact that $\alpha=\fr{\log\rho}{\log\lambda}-1$.

Invariance by $g$ is easily checked.$ \Box $

\ \\
Thus, we have a quasiconformal structure on a symmetric open set
 (image of a disk $D_R$ under $\ph^{ -1}$),
  which is invariant by $g$
 with symmetric complex dilatation coefficient
 $\mu_{\rho}$  
 (i.e. $\mu_{\rho}(1/\barre{z}) = \barre{\mu_{\rho}}(z)$).

We propagate $\mu_{\rho}$ over almost all the attracting basin by pull back by $g$.
The resulting Beltrami form is still symmetric and we can extend it by
 $0$
 outside the basin.

To do so, the Beltrami form $\sigma_{\rho}=\mu_{\rho}\fr{d\barre{z}}{dz}\in L ^\infty$
 is constructed in this manner (compare \cite{ShishiQC}) :
\begin{itemize}
\item $\sigma_{\rho}=\sigma_{std}$, the standard conformal structure of
 $\mathbb{P}^1$, outside the Borel set :
\rae
 \{z/\forall n\in\naturels,\,\forall m\in\naturels,\,g^{\circ n}(z)\nin g^{\circ m}(U)\},
\are
\item on $U$, $\sigma_{\rho}$ is the one given by lemma
 \ref{lem sur le changement de structure c du disque num2}, it si invariant by $g$,
\item on the preimages $g^{-n}(U)$ of $U$, we just pull $\sigma_{\rho|U}$ back by $g$,
 $\mu_{\rho}(z)=
\fr
{\barre{g^{\circ n}{}'(z)}}
{g^{\circ n}{}'(z)}
\mu_{\rho}(g^{\circ n}(z))$
 where $n$ is such that
 $g^{\circ n}(z)\in U$.
\end{itemize}
Since $U$ is included in the immediate basin of an attracting cylce and
 $\sigma_{\rho|U}$ is invariant by $g$,
 $\sigma_{\rho|U}$ 
 is uniquely determined by any restriction on any fundamental annulus of the dynamics of $g$ on $U$.

Let
 $\sigma_{\rho} =
 \mu_{\rho}\fr{d\barre{z}}{dz}$.
From holomorphy of $g$ it follows that
$||\mu_{\rho}||_\infty=
||\mu_{\rho}||_{L^\infty(\mathbb{P}^1)}
=||\mu_{\rho}||_{L^\infty(U)}$.

The coefficient $\mu_{\rho}$ (or $\sigma_{\rho}$)
 is easily seen to be symmetric on $U$. Moreover,
 $\forall n$, $\forall z$, 
 $g^{\circ n}{}'(1/\barre{z})=\barre{g^{\circ n}{}'(z)}$.
 It follows that for $z$ such that $g^{\circ n}(z)\in U$,
\rae
\mu_{\rho}(1/\barre{z})=
\fr{\barre{g^{\circ n}{}'(1/\barre{z})}}{g^{\circ n}{}'(1/\barre{z})}\mu_{\rho}(g^{\circ n}(1/\barre{z}))
=\barre{\left(\fr{\barre{g^{\circ n}{}'(z)}}{g^{\circ n}{}'(z)}\mu_{\rho}(g^{\circ n}(z))\right)} 
 =  \barre{\mu_{\rho}(z)}.
\are
\ \\

\begin{rem}
The above pullback discussion is independent of $\rho$.
 As a consequence, the Beltrami form $\sigma_{\rho}$ depends analytically on $\rho$.
\end{rem}

Since
 $||\mu_\chi||_\infty=||\mu_{\rho}||_\infty <1$,
 the Riemann measureable mapping theorem yields the existence of a unique quasiconformal homeomorphism
 $\Phi_\rho:(\mathbb{P}^1,0,\infty,x)\longrightarrow(\mathbb{P}^1,0,\infty,x)$
 with dilatation $\mu_{\rho}$.
Moreover, this homeomorphism is symmetric with respect to the circle (compare section
 \ref{sec gen prop of gab})
 and by the regularity of the dependence of the solutions of the Beltrami equation with respect to the parameters theorem
 (compare for example \cite{AhlforsQCbook}),
 $\rho\mapsto\Phi_\rho(z)$
 is real-analytic for any $z\in\mathbb{P}^1$.

From the dilatation of $\Phi_{\rho}$, it follows that 
 $\tilde{g} \dpe \Phi_{\rho} \circ g \circ \Phi_{\rho}^{ -1}$ and 
 $\tilde{\ph}\dpe\left(\Phi_{\rho}\circ\ph^{-1}\circ\chi\right)^{-1}$
 are holomorphic on their respective domains of definition :
\begin{equation*}
\tilde{\ph} : \Phi_{\rho}(U) \stackrel{\Phi_{\rho}^{-1}}{\longrightarrow}U \stackrel{\ph}{\longrightarrow} D_R
 \stackrel{\chi}{\longrightarrow}D_r
\end{equation*}

\begin{equation}
\xymatrix{
\Phi_{\rho}(U) \ar[r]^{\Phi_{\rho}^{-1}} \ar[d]_{\tilde{g}}  &
 U \ar[d]_g \ar[r]^{\ph}  
& D_R \ar[r]^{\chi} \ar[d]^{\times\lambda}  & D_r \ar[d]^{\times\rho}\\
\Phi_{\rho}(h(U)) \ar[r]^{\Phi_{\rho}^{-1}} & h(U)\ar[r]^{\ph}&\lambda D_R \ar[r]^{\chi} & \rho D_r
}
\end{equation}

The new mapping $\tilde{g}$ has an attracting cycle of period $p$, like $g$,
 with a local linearizing mapping
 $\tilde{\ph}$ defined in the neighbourhood of the point 
 $\Phi_{\rho}(x)$ (which belongs to this cycle). The new multiplier is thus $\rho = \lambda^{1+\alpha}$.
Note that all values in the interval
 $]0,1[$
 can be assigned to
 $\rho$.

We need to see that we obtained a mapping which belongs to the family of complexified double standard maps
 $(g_{a,b})_{a,b}$.
 In the next section, we will see that, after conjugating by a suitable rotation, the new map belongs to the family.

\subsection{Back to the double standard family}
We now check that the mapping $\tilde{g}$
 constructed in the previous section, is conjugated by a rotation
 $R$
 to a unique element
 $g_{\tilde{a},\tilde{b}}$
 in the family of mappings $(g_{a,b})_{a,b}$.
This is a consequence of the following and of the discussion in section \ref{sec gen prop of gab}
about the existence and unicity of the conjugacy by a rotation.

\begin{prop}\label{lem le resultat est dans la famille important de tourner}
Let $(a,b,c)\in\complexes\times\complexes^*\times\complexes^*$ and let
 $g_{a,b,c}:\complexes^*\rightarrow\complexes^*$ be defined by :
\rae
g_{a,b,c}(z) & = & e^{2\ii\pi a}z^2 e^{-\left(bz-c/z\right)}.
\are
Let $\ph$, $\psi$ be orientation preserving homeomorphisms of the Riemann sphere $\projectif^1$
 fixing $0$ and $\infty$.

If $\psi\circ g_{a,b,c} \circ\ph:\complexes^*\rightarrow\complexes^*$ is holomorphic, then there exist
 $(\alpha,\beta,\gamma)\in\complexes\times\complexes^*\times\complexes^*$\linebreak
 such that
 $\psi\circ g_{a,b,c} \circ\ph=g_{\alpha,\beta,\gamma}$.

If we suppose that $g_{a,b,c}$ is a symmetric map with respect to the circle, i.e.
 $g_{a,b,c}\left(1/\barre{z}\right)=1/\barre{g_{a,b,c}(z)}$,
 hence that $a\in\reels$ and $c=\barre{b}$, and that $\ph$ and $\psi$
 are also symmetric with respect to the circle, then the mapping $g_{\alpha,\beta,\gamma}$ is also symmetric,
 $\alpha\in\reels$, $\gamma=\barre{\beta}$, and we have $|\beta|<1$ if and only if $|b|<1$.
\end{prop}

\preuve
Without loss of generality we can suppose $\ph$ and $\psi$ are quasiconformal.

Indeed, let $V$ be the set of critical values of $g_{a,b,c}$.
It contains one or two point-s and $\psi(V)$ is the set of critical values of $h\dpe\psi\circ g_{a,b,c}\circ \ph$.
From the finiteness of $V$ and the compactness of $\projectif^1$, it follows that there exists an isotopy
 $\psi_t:\projectif^1\rightarrow\projectif^1$ between $\psi=\psi_0$ and a quasiconformal homeomorphism
 $\psi_1$ such that, for all $t\in[0,1]$, $\psi_t(0)=0$, $\psi_t(\infty)=\infty$ and such that for all $v\in V$,
 $\psi_t(v)=\psi(v)$.

Thanks to this last property, the isotopy $\psi_t$ can be lifted up into an isotopy $\ph_t$ starting from
 $\ph_0=\ph$. The final function $\ph_1$
 of this isotopy is quasiconformal since it is locally a composition of conformal or quasiconformal mappings.

Since $\ph$ and $\psi$ are orientation preserving homeomorphisms of $\complexes^*$
 they induce the identity map on the first homology group of $\complexes^*$.
 This group is generated by the class of any curve winding one time round $0$ in the anticlockwise direction.
 As a consequence
\rae
\inv{2\ii\pi}\int_{(0,+)}\fr{h'(z)}{h(z)}dz & =& \inv{2\ii\pi} \int_{(0,+)}\fr{g_{a,b,c}'(z)}{g_{a,b,c}(z)}dz.
\are
The latter is equal to $2$, from what follows there exists a holomorphic function
 $u:\complexes^*\rightarrow\complexes^*$
 such that $h(z)=z^2 e^{u(z)}$.
In order to characterize the function $u$, we need to know the growth order of $h$ near $0$ and $\infty$.

From now on, we suppose that $\ph$ and $\psi$ are $K$-quasiconformal homeomorphisms of $\projectif^1$.
It is well known that such mappings are locally Hölder continuous of exponent $1/K$.
So, by applying this to the inverse mapping of $\ph$ near $\infty$, we can find $R_\infty>0$ and
 $C_\infty>0$ such that for all $z$ satisfying $|\ph(z)|\geq R_\infty$, we have $|\ph(z)|\leq C_\infty|z|^K$.
Moreover, we can fix $R_\infty$ and $C_\infty$ such that the same holds if we only change $\ph$
 and replace it by $\psi$.

We concentrate on estimating the growth order of $h$ near $\infty$. Since $\ph(z)\tend\infty$ when
 $z\tend\infty$, we can suppose that $|z|$ is such that $|\ph(z)|\geq R_\infty$.
 Then
\rae
 \parel\left(-\left(b\ph(z)+\fr{c}{\ph(z)}\right)\right) & \leq & |b| |\ph(z)| + \fr{|c|}{|\ph(z)|}\\
 & \leq & |b|C_\infty |z|^K + \fr{|c|}{R_\infty}.
\are
Hence $|g\circ\ph(z)|\leq C|z|^{2K}e^{C'|z|^K}$, for some positive constants $C$, $C'$.
Therefore if
$ |\psi\circ g_{a,b,c}\circ \ph(z)| \geq R_\infty$, we have
$ |h(z)| = |\psi\circ g_{a,b,c}\circ \ph(z)| \leq C'' |z|^{2K^2} e^{C''' |z|^K}$ for some other positive constants
$C''$, $C'''$. 
Thus the growth order of $h$ near $\infty$ is finite.

This conclusion also holds near $0$ from a similar argument using Hölder continuity of $\ph$ near $0$.

Therefore the Laurent series development of $u$
 is a non trivial sum of a finite number of relative powers of $z$,
\rae
u(z) & = & \sum_{n=-p}^q a_n z^n,
\are
with $(p,q)\in\naturels^*\times\naturels^*$, $a_{-p}\neq 0$ and $a_q\neq 0$.

From the definition of $h$, it follows that $h$ has the same number of critical points as $g_{a,b,c}$.
From section \ref{sec gen prop of gab}\linebreak we know that the mapping $g_{a,b,c}$
 has two critical points counting multiplicity.

We have
\rae
h'(z) & = & \left(2z+ z^2u'(z)\right)e^{u(z)} \\
 & = & \inv{z^{p-1}}\left(2z^p+z^{p+1}u'(z)\right)e^{u(z)}.
\are
The function $P(z)=\left(2z^p+z^{p+1}u'(z)\right)$ is a polynomial of exact degree $p+q$ such that
 $P(0)\neq 0$, thus $h$ has $p+q$ critical points counting multiplicity. Therefore we have $p=q=1$.
This ends the proof of the first part of proposition
 \ref{lem le resultat est dans la famille important de tourner}.

The symmetry part of the last assertion is straightforward, while the equivalence on the bound $|\beta|<1$
 follows from the discussion in the beginning of section \ref{sec gen prop of gab}$\square$

\subsection{Type of the deformed cycle}

Given any $\rho\in]0,1[$,
 we have a process for constructing a map in the double standard family by deformation of a map belonging to a tongue, such that the new map has an attracting cycle with multiplier
 $\rho$.
We check that this construction does not change the type.

\begin{prop}
The tuple of parameters $(\tilde{a},\tilde{b})$
 corresponding to the new function $\tilde{g}$ is of type $\tau$.
\end{prop}

\preuve

Let's recall that $\tau$ is the image  by the map
 $\phi_{\tilde{a},\tilde{b}}(x)=\ds{\lim_{n\tend \infty}}F_{a,b}^{\circ n}(x)/2^n$ (see \cite{Misurod1})
 of the point $\tilde{x}$
 of the attracting cycle of
 $g_{\tilde{a},\tilde{b}}$
 which belongs to the basin component that contains the critical points.

The proposition is a consequence of the following uniqueness property of the semiconjugacy $\phi_{a,b}$.

\begin{lem}\label{unicite de la semiconjuguaison}
Let $f:\tore^1\rightarrow \tore^1$ be a monotonic continuous mapping.
 If $\ph,\psi:\tore^1\rightarrow \tore^1$ are non decreasing continuous mappings of degree
 $1$
 such that 
 $\ph \circ f=2\times\ph$ and $\psi \circ f=2\times\psi$,
 then $\ph=\psi$.
\end{lem}

\begin{rem}
In the above lemma, the continuity assumption is not necessary, but, if one considers discontinuous mappings, the  "non decreasing" assumption has to be restated in order to take that specific context into account.
\end{rem}

\preuve

Let $F$ be a lift of $f$ and let $\tilde{\ph}$ be a lift of $\ph$.
Since $\tilde{\ph}\circ F$ is an non decreasing real function,
 it is a lift of $2\ph$.

Therefore, there exists an integer constant $k$ such that 
\rae
\forall x\in\reels,\,\tilde{\ph}(F(x))=2\tilde{\ph}(x)+k.
\are
It follows that if we take $\ph_1=\tilde{\ph}+k$, then we have
 $\ph_1\circ F=2\ph_1$.
We assume we have the same properties for a lift $\psi_1$ of $\psi$.

The mappings $\ph_1$ and $\psi_1$  are both non decreasing of degree $1$
 and locally bounded, thus, the function $\ph_1 - \psi_1$
 is periodic and bounded. But then
\rae 
\ph_1(x)-\psi_1(x)& =& \demi\left(\ph_1(F(x))-\psi_1(F(x))\right),\\
\are
so
\rae
\ph_1(x)-\psi_1(x)&=& \left(\demi\right)^n\left(\ph_1(F^{\circ n}(x))-\psi_1(F^{\circ n}(x))\right),
\are
for all $n$. Consequently $\ph_1(x)-\psi_1(x)=0$.
 $\Box$

%
%
\ \\

The following diagram is commutative 
\begin{equation}
\xymatrix{
\tore^1\ar[d]_{f_{\tilde{a},\tilde{b}}}\ar[r]^{\Psi^{-1}}
 & \tore^1\ar[d]_{f_{a,b}}\ar[r]^{\phi_{a,b}} & \tore^1 \ar[d]_D \\
\tore^1\ar[r]_{\Psi^{-1}} & \tore^1\ar[r]_{\phi_{a,b}} & \tore^1 }
\end{equation}
where $\Psi$ is defined by $\exp\circ\Psi = \Phi_{|\sphere^1}\circ\exp$.
Therefore, we have $\phi_{\tilde{a},\tilde{b}}=\phi_{a,b}\circ\Psi^{-1}$.
As a consequence $(\tilde{a},\tilde{b})$ is of type $\tau$.

\section{Path}
\label{section chemin dans la langue que c'est un chemain}

Lemma \ref{lem sur le changement de structure c du disque} gives the value of the multiplier
 $\rho$
 of the new mapping
 $g_{\tilde{a},\tilde{b}}$, it is $\rho=\lambda^{1+\alpha}$,
 where
 $\alpha$ is taken in the interval $ ]-1,+\infty[$.
 It is then clear that all values in the interval
 $]0,1[$
 can be assigned to
 $\rho$
when we change the parameter of our construction
 (it suffices to set $\alpha=\fr{\log\rho}{\log\lambda}-1$).

Define
 $\gamma(\rho)\dpe(\tilde{a}(\rho),\tilde{b}(\rho))$
 where $(\tilde{a}(\rho),\tilde{b}(\rho))$
 is the tuple of parameters corresponding to
 $f_{\tilde{a},\tilde{b}}$, which is the restriction of $g_{\tilde{a},\tilde {b}}$
 to the unit circle (where
 $g_{\tilde{a},\tilde {b}}$
 was constructed from 
 $g_{a,b}$
 by the process of quasiconformal deformation decribed above).

This mapping
 $\gamma$,
 is well defined on $]0,1[$,
 satisfies
 $\gamma(\lambda)=(a,b)$ 
 and takes values into the type $\tau$ tongue $T_\tau$.

This mapping $\rho\mapsto \gamma(\rho)$
 will be shown to be continuous. Indeed it is real analytic.
Afterwards, it remains to show that $\gamma(0)$
can be defined and corresponds to a parameter with a superattracting cycle of type
$\tau$.

\subsection{Analyticity}\label{subsection analycite du chemin}
The proof of the analytic dependence in $\rho$
 of $(\tilde{a},\tilde{b})=(a(\rho),b(\rho))$
 is based on the study of the regularity of the dependence of the critical points and critical values of
 $\tilde{g}$,
 the resulting map of the quasiconformal deformation, and of the rotation which turns
 $\tilde{g}$
 into a member of the family
 (cf. lemma \ref{lem le resultat est dans la famille important de tourner}).

Recall that $\rho\mapsto \Phi_{\rho}(z)$ is a real-analytic map for all $z$.
The critical points of
 $\tilde{g}_\rho=\Phi_\rho\circ g\circ \Phi_{\rho}^{-1}$, 
 which are topologically characterized, are the images of the critical points of
 $g$
 by
 $\Phi_\rho$.
Let
 $\omega$
 be one of the two critical points of
 $g$
 and
 $\tilde{\omega}_\rho\dpe\Phi_\rho(\omega)$.
Since
 $\rho\mapsto\Phi_\rho(\omega)$
 is real-analytic, the critical point 
 $\tilde{\omega}_\rho$
 depends analytically on
 $\rho$.

The aim of the conjugacy by the rotation
 $R_\rho$
 is to normalize
 $\tilde{g}_\rho$
 so as to have a mapping
 $g_\rho=R_\rho\circ \tilde{g}_\rho\circ R_\rho^{-1}$
 belonging to the complexified double standard family.
Since the critical points of
 $\tilde{g}_\rho(z)=\lambda z^2 e^{-(bz-\barre{b}/z)}$
 are
 $\fr{1\pm\sqrt{1-|b|^2}}{b}\in\inv{b}\reels$,
 it comes to move the critical point on the real positive half line.
Specifically,
 $R_\rho z = \fr{|\tilde{\omega}_\rho|}{\tilde{\omega}_\rho}z$.

Therefore, the corresponding critical point of
 $g_\rho$ : 
 $\omega_\rho=R_\rho\tilde{\omega}_\rho=|\tilde{\omega}_\rho|$,
 is still depending analytically on
 $\rho$.

From the formulas of the critical points of
 $g_\rho=g_{a(\rho),b(\rho)}$,
 we get
\rae
 b(\rho)=\fr{2\omega_\rho}{1+\omega_\rho^2},
\are 
 which is true for both of the critical points $\omega_\rho$.

We can express
 $a(\rho)$
 in terms of
 $\rho$
 by using the expression of the critical value of
 $g_\rho$
 which is at the same time image of
 $\omega_\rho$
 by
 $g_\rho(z)=e^{2\pi a(\rho)\ii}z^2 e^{-b(\rho)(z-1/z)}$
 and image of
 $h(\omega)$
 by
 $R_\rho\circ\Phi_\rho$. We have
\rae
 e^{2\pi a(\rho)\ii}=
\fr{R_\rho\circ\Phi_\rho(f(\omega))}{\omega_\rho^2 e^{-b(\rho)(\omega_\rho-1/\omega_\rho)}}.
\are

 So,
 $\rho\mapsto(a(\rho),b(\rho))$
 is an analytic path from
 $]0,1[$
 into
 $T_\tau\subset \tore^1\times[0,1]$.

\subsection{Ending of the path when the multiplier tends to $0$}
\label{section, aboutissement du chemin}

We have shown there exists a continuous path
 $\rho\in]0,1[\mapsto \gamma(\rho)=(a(\rho),b(\rho))$
 in the type
 $\tau$ tongue,
 along one direction of which the multiplier of the type
 $\tau$ cycle is $\rho$ and tends to zero.
We have to show that this path ends up at a defined limit when $\rho\tend 0$,
 that this limit corresponds to a type
 $\tau$
 parameter with superattracting cycle and that it doesn't depend on the starting point in the tongue.

\ \\

We first show the limit
 $\ds{\lim_{\rho\tend 0}}(a(\rho),b(\rho))$
 exists.

Note that $b(\rho)\leq 1$ for all $\rho$ since for $b>1$, $f_{a,b}$
 is not order preserving and thus cannot be topologically conjugated to an order preserving map.

Since
 $f_{a,b}'(x)\geq 2(1-b)$
 for all
 $(a,b,x)$
 and
 $(f_{a(\rho),b(\rho)}^{\circ p})'(x(\rho))=\rho\ds{\tend_{\rho\tend 0}}0$,
 where
 $x(\rho)$
 is any of the points of the attracting cycle (of period
 $p$)
 of
 $f_{a(\rho),b(\rho)}$,
The value of
 $b(\rho)$
 tends to
 $1$
 when
 $\rho\tend 0$.
In paticular, from the continuous dependence of the semiconjugacy
 $\phi_{a,b}$ with respect to $(a,b)$,
 it follows that for any limit point of
 $(a(\rho),b(\rho))$
 when
 $\rho\tend 0$
 the corresponding map has a superattracting cycle of exact period
 $p$.

The set
 $\Lambda$
 of limit points of
 $\gamma(\rho)=(a(\rho),b(\rho))$
 when
 $\rho\tend 0$
 is then included in the set of points
 $(a,1)\in\tore^1\times\{1\}$.
This set
 $\Lambda$
 is compact and connected since it is the decreasing intersection of the connected compact sets
$\barre{\gamma(]0,1/n[)}$.

Thus
 $\Lambda$
 is a point or contains an open set
 $]a_1,a_2[$
 in
 $\tore^1\times\{1\}$. We will prouve that its interior is empty.

To poceed we need the following lemma, also used later.
\begin{lem}\label{lem, tout petit sur monotonie etc.}
For any
 $p\in\naturels$, 
 the mapping
 $\tore^1\ni a\mapsto f_{a,1}^{\circ p}(1/2)\in\tore^1$
 is strictly increasing and of degree
 $2^p-1$.
\end{lem}

\preuve
Monotonicity comes from the fact that
 $\fr{\partial}{\partial a} \left(F_{a,1}^{\circ p}(1/2)\right)\geq 1$,
 which can be seen by a direct calculation (compare \cite{Misurod2}).

Now, it is sufficent to see that for all $(a,b,x)\in\reels\times[0,1]\times\reels$ and all
 $p\in\naturels$ :
\rae
 F_{a+1,b}^{\circ p}(x)=F_{a,b}^{\circ p}(x)+2^p-1
\are
 (recall that
 $f_{a,b}:\tore^1\rightarrow\tore^1$
 is the quotient in
 $\tore^1$
 of the function
 $F_{a,b}:\reels\rightarrow\reels$).

This is true for $p=1$, moreover, an induction argument on
 $p$ gives :
\rae
 F_{a+1,b}^{\circ p+1}(x) & = & F_{a+1,b}(F_{a+1,b}^{\circ p}(x))\\
 & = & F_{a,b}(F_{a,b}^{\circ p}(x)+2^p-1)+1\\
 &= & F_{a,b}^{\circ p+1}(x) + 2^{p+1}-1.\Box
\are

Now we see that the limit is unique for otherwise,
 $1/2$, which is the only critical point of $g_{a,1}$,
 would be a $p$-periodic point of
 for all of the maps $f_{a,1}$ such that
 $a\in]a_1,a_2[$,
which is impossible by the above.
The limit is thus unique.

\ \\

Secondly, the type of the limit $(a,1)$ must be
 $\tau$.
Indeed it has a type $\tau'$. But since the set of parameters of a given type is an open set in
 $\tore^1\times[0,1]$, all the values of the path close enough to $(a,1)$
 have the same type $\tau'$ which implies that $\tau'=\tau$.

\ \\

Finally, we have a type 
 $\tau$
 parameter with a superattracting cycle, and this parameter is the limit point of the path 
 $\gamma(\rho)$
 when
 $\rho\tend 0$.
Now we have to show that this parameter is unique.

We know from lemma \ref{lem, tout petit sur monotonie etc.}
 that for any given period
 $p\in\naturels$,
 there exist exactly
 $2^p-1$
 parameters
 $(a,1)\in\tore^1$
 for which the mapping
 $f_{a,b}$
 has a superattracting cycle of period dividing
 $p$.

Moreover, for each type
 $\tau\in\tore^1$
 (periodic point of $D:x\mapsto 2x \mod 1$),
 the type $\tau$ tongue
 $T_\tau$ 
 is nonempty, it comes from the fact that
 $\tore^1 \ni a \mapsto \phi_{a,1}(1/2)\in\tore^1$
 is monotonic of degree $1$,
 which implies that the critical point belong to the attracting basin of any type of cycle.

Therefore, for all such
 $\tau$,
 there exists at least one parameter
 $(a,1)\in\tore^1\times\{1\}$
 of type
 $\tau$
 with superattracting cycle.
Since
 $D$
 has exactly
 $2^p-1$
 periodic points of period dividing
 $p$
 (namely, the numbers
 $\fr{k}{2^p-1} \mod1$ with
 $k=0,\dots,2^p-1$),
 the type
 $\tau$
 parameters with superattracting cycle are unique for all
 $\tau$.

\section{Open questions}


The other question left unaswered in the reference \cite{Misurod2}
 was about the order of contact of the left and right boundaries of the tongue.
  The order of contact might be
 $\demi$
 for all of the tongues.
 See this reference for a study of examples.
 We have reasons to believe that the order of contact is the same for all of the tongues.


When we consider the dynamics of the family
 $g_{a,b}$ instead of $f_{a,b}$,
 we see that there can also be a symmetric pair of attracting cycles in $\complexes^*$.
We might expect the corresponding parameters to form Mandelbrot-like families or even tricorn shapes in the parameter space $\tore^1\times[0,1]$. 
Mandelbrot like shapes can be seen (see pictures below) attached on each tongue, with similar non-local connectedness as previous examples (compare \cite{NakaneSchleicher}, \cite{EpsteinYampolski}).
We should notice that the complexified double standard family may be related to a rational Blashke product with same kind of dynamics near the circle (the restriction to the circle is of degree two and strictly monotonic), see \cite{FagellaGarijo}.

\begin{center}
\begin{figure}
\includegraphics[scale=0.04]{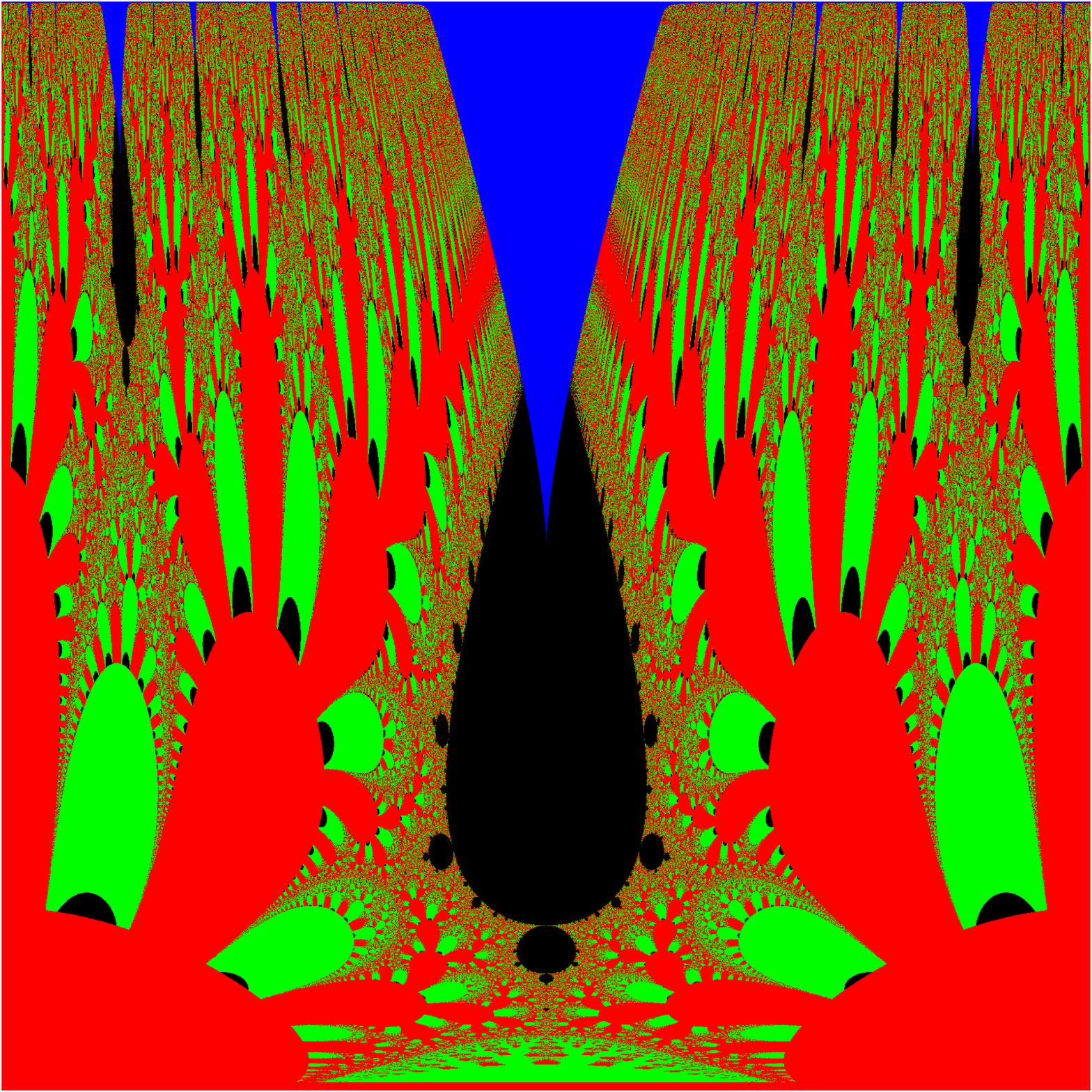}
\includegraphics[scale=0.106]{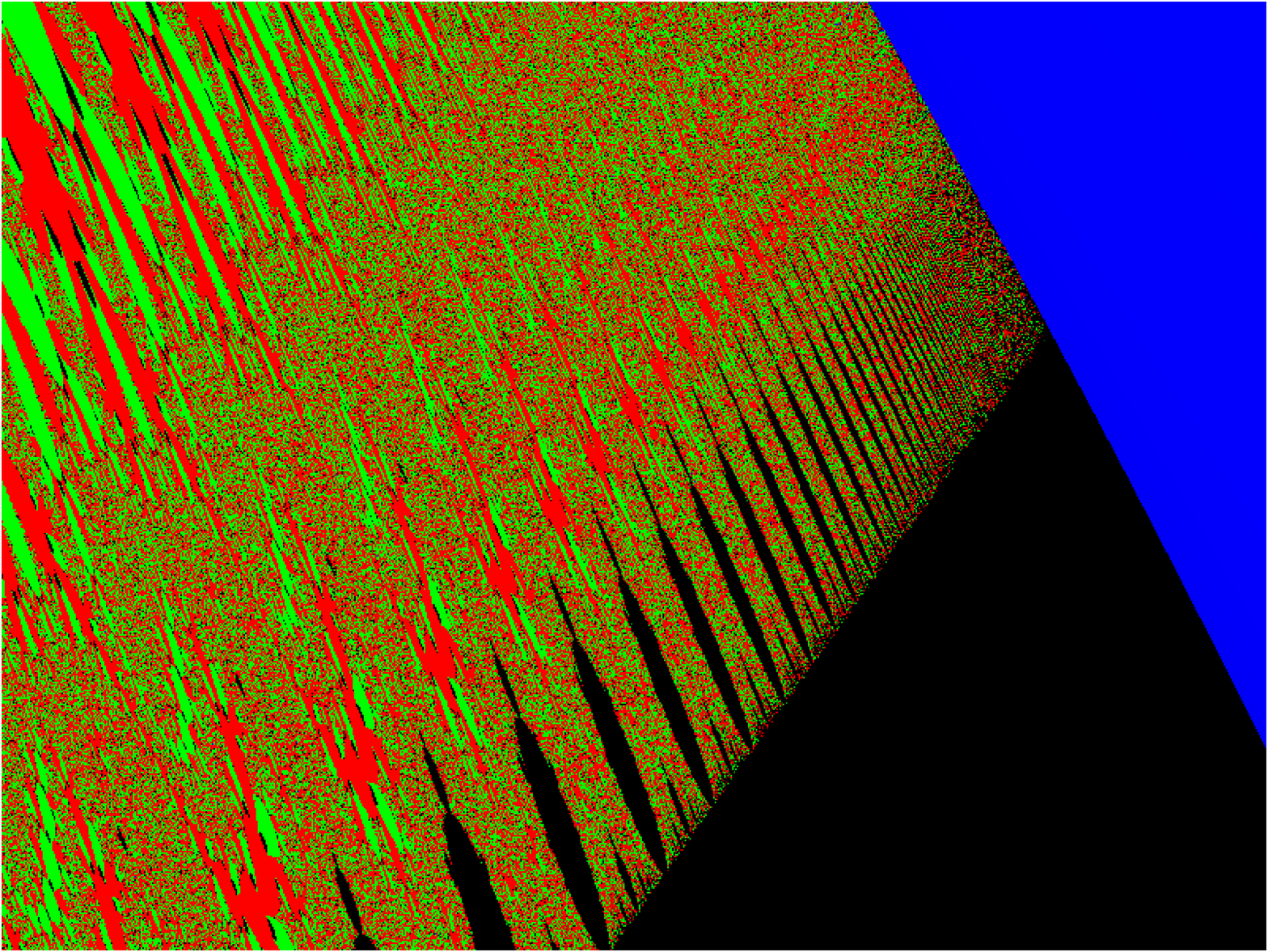}
\caption{Parameter plane for the family $(g_{a,b})_{a,b}$, tongues are blue, parameter with escaping critical points are red and green, and the others,  the parameters for which the critical point stays in a neighbourhood of the cricle without tending to it after some iterations, are black. On right a magnification of the contact between the Mandelbrot like shape and the tongue.}
\end{figure} 
\end{center}

\bibliographystyle{plain}
\bibliography{mabib}

\end{document}